\documentclass[a4paper]{llncs}
\usepackage{graphicx}
\graphicspath{{figures/},{figures/fig1/}}
\usepackage{custom_style}  
\usepackage{bm}
\usepackage{amsmath}
\newtheorem{df}{Definition}
\newtheorem{theo}{Theorem}
\newtheorem{rk}{Remark}
\usepackage{amsfonts}
\usepackage{amssymb}
\usepackage{graphicx}
\usepackage{epstopdf}
\usepackage{tikz}
\usepackage{algorithmic}
\usepackage{multirow}
\usepackage{url}
\usepackage{relsize}
\usepackage{paralist}
\usepackage{float}
\usepackage{subcaption}
\usepackage[ruled,boxruled,nofillcomment]{algorithm2e}
\usepackage{cleveref}
\usepackage{comment}
\def\*#1{\mathbf{#1}}

\begin{document}

\title{ECKO: Ensemble of Clustered Knockoffs for robust multivariate
  inference on fMRI data}
%
%
%
\author{Tuan-Binh Nguyen \inst{1,2,3} \and
J\'er\^ome-Alexis Chevalier \inst{1,2,4} \and
Bertrand Thirion\inst{1,2}}
\authorrunning{F. Author et al.}
%
\institute{
  Parietal team, Inria Saclay, France \and
  CEA/Neurospin, Gif-Sur-Yvette, France \and
  LMO - Laboratoire de Math\'ematiques d'Orsay \and
  Telecom Paristech, Paris, France
}
%
\maketitle              

\begin{abstract}

Continuous improvement in medical imaging techniques allows the acquisition of
higher-resolution images.
When these are used in a predictive setting, a greater number of explanatory
variables are potentially related to the dependent variable (the response).
Meanwhile, the number of acquisitions per experiment remains limited.
In such high dimension/small sample size setting, it is desirable to find the
explanatory variables that are truly related to the response while controlling
the rate of false discoveries.
To achieve this goal, novel multivariate inference procedures, such as knockoff
inference, have been proposed recently.
However, they require the feature covariance to be well-defined, which is
impossible in high-dimensional settings.
In this paper, we propose a new algorithm, called Ensemble of Clustered
Knockoffs, that allows to select explanatory variables while controlling the
false discovery rate (FDR), up to a prescribed spatial tolerance.
The core idea is that knockoff-based inference can be applied on
groups (clusters) of voxels, which drastically reduces the problem's
dimension; an ensembling step then removes the dependence on a fixed
clustering and stabilizes the results.
We benchmark this algorithm and other FDR-controlling methods on brain imaging
datasets and observe empirical gains in sensitivity, while the false
discovery rate is controlled at the nominal level.
\end{abstract}

\section{Introduction}  
\label{sec:intro}
Medical images are increasingly used in predictive settings, in which
one wants to classify patients into disease categories or predict some
outcomes of interest.
Besides predictive accuracy, a fundamental question is that of
\textit{opening the black box}, \ie understanding the combinations of
observations that explains the outcome.
%
%
%
A particular relevant question is that of the importance of image
features in the prediction of an outcome of interest, conditioned on
other features.
Such conditional analysis is a fundamental step to allow causal inference on
the implications of the signals from image regions in this outcome; see e.g.
\cite{weichwald_causal_2015} for the case of brain imaging.
However, the typical setting in medical imaging is that of high-dimensional
small-sample problems, in which the number of samples \( n \) is
much smaller than the number of covariates \( p \).
This is further aggravated by the steady improvements in data resolution.
In such cases, classical inference tools fail, both theoretically and
practically.
One solution to this problem is to reduce the massive number of covariates by
utilizing dimension reduction, such as clustering-based image
compression, to reduce the number of features to a value close to n; see e.g.
\cite{buhlmann_correlated_2012}.
%
%
This approach can be viewed as the bias/variance
trade-off: some loss in the localization of the predictive features
---bias--- is tolerated as it comes with less variance ---hence higher
power--- in the statistical model.
This is particularly relevant in medical imaging, where localizing
predictive features at the voxel level is rarely important: one is
typically more interested in the enclosing region.\\
\indent However, such a method suffers from the arbitrariness of the
clustering step and the ensuing high-variance in inference results
with different clustering runs, as shown empirically in
\cite{chevalier_ensembling_2018}.
\cite{chevalier_ensembling_2018} also introduced an algorithm called
Ensemble of Clustered Desparsified Lasso (ECDL), based on the
inference technique developed in \cite{zhang_confidence_2014}, that
provides p-values for each feature, and controls the Family Wise Error
Rate (FWER), i.e. the probability of making one or more false
discoveries.
In applications, it is however more relevant to control the False Discovery Rate
(FDR) \cite{benjamini_controlling_1995}, which indicates the expected fraction
of false discoveries among all discoveries, since it allows to detect a greater
number of variables.
In univariate settings, the FDR is easily controlled by the Benjamini-Hochberg
procedure \cite{benjamini_controlling_1995}, valid under independence or
positive correlation between features.
It is unclear whether this can be applied to multivariate statistical settings.
A promising method which controls the FDR in multivariate settings is the
so-called knockoff inference \cite{barber_controlling_2015,candes_panning_2016},
which has been successfully applied in settings where $n \approx p$.
However, the method relies on randomly constructed knockoff
variables, therefore it also suffers from instability.
Our contribution is a new algorithm, called Ensemble of Clustered
Knockoffs (ECKO), that \textit{i)} stabilizes knockoff inference
through an aggregation approach; \textit{ii)} adapts knockoffs to $n
\ll p$ settings.
This is achieved by running the knockoff inference on the reduced
data and ensembling the ensuing results.\\
\indent The remainder of our paper is organized as follows: \Cref{sec:theory}
establishes a rigorous theoretical framework for the ECKO algorithm;
\Cref{sec:experiment} describes the setup of our experiments with both synthetic
and brain imaging data predictive problems, to illustrate the performance of
ECKO, followed by details of the experimental results in \Cref{sec:results};
specifically, we benchmark this approach against the procedure proposed in
\cite{gaonkar_interpreting_2015}, that does not require the clustering step, yet
only provides asymptotic ($n\rightarrow\infty$) guarantees.
We show the benefit of the ECKO approach in terms of both statistical
control and statistical power.
%
\section{Theory}  
\label{sec:theory}
%
\subsection{Generalized Linear Models and High Dimensional Setting}
\label{sec:linear_model}
%
Given a design matrix $\*X \in \bbR^{n \times p}$ and a response vector
$\*y \in \bbR^n$, we consider that the true underlying model is of the following
form:
\begin{equation} \label{linear_model}
  \*y = f(\*X \*w^*) + \sigma \bm\epsilon \enspace ,
\end{equation}
where $ \*w^* \in \bbR^p $ is the true parameter vector, $ \sigma \in \bbR^+ $
the (unknown) noise magnitude, $\bm\epsilon \sim \cN(\mathbf{0}, \bI_n) $ the
noise vector and $f$ is a function that depends on the experimental setting
(e.g. $f = Id$ for the regression problem or e.g. $f=sign$ for the
classification problem).
The columns of $\*X$ refer to the explanatory variables also called features,
while the rows of $\*X$ represent the coordinates of different samples in the
feature space.
We focus on experimental settings in which the number of features $p$ is much
greater than the number of samples $n$ \ie $p \gg n$.
Additionally, the (true) support denoted by $S$ is given by $S = \{ k \in [p] : \*w_k^* \neq 0 \}$.
Let $\hat{S}$ denotes an estimate of the
support given a particular inference procedure.
We also define the signal-to-noise ratio (SNR) which allows to assess the noise
regime of a given experiment:
\begin{equation}\label{eq:snr}
  \text{SNR} = \dfrac{\norm{\*X \*w^* } ^2_2}{\sigma^2 \norm{\bm\epsilon}^2_2} \enspace .
\end{equation}
A high SNR means the signal magnitude is strong compared to the noise, hence it
refers to an easier inference problem.
%
\subsection{Structured Data}
\label{sec:structured data}
%
In medical imaging and many other experimental settings, the data stored in the
design matrix $\*X$ relate to \emph{structured} signals.
More precisely, the features have a peculiar dependence structure that
is related to an underlying spatial organization, for instance the
spatial neighborhood in 3D images.
Then, the features are generated from a random process acting on this underlying
metric space.
In our paper, the distance between the $j$-th and the $k$-th features is denoted
by $d(j,k)$.
%
\subsection{FDR control}
\label{sec:FDR}
%
In this section, we introduce the false discovery rate (FDR) and a spatial
generalization of the FDR that we called $\delta$-FDR.
This quantity is important since a desirable property of an inference procedure
is to control the FDR or the $\delta$-FDR.
In the following, we assume that the true model is the one defined in \Cref{sec:linear_model}.
%
%
\begin{df}\label{df:FDP} \emph{False discovery proportion (FDP).}
  Given an estimate of the support $\hat{S}$ obtained from a particular
  inference procedure, the false discovery proportion is the ratio of the number
  selected features that do not belong to the support (false discoveries)
  divided by the number of selected features (discoveries):
\begin{equation}\label{eq:fdp}
\text{FDP} = \frac{\#\{k \in \hat{S}: k \notin S\}}{\#\{k \in \hat{S}\}}
\end{equation}
\end{df}
\begin{df}\label{df:delta_FDP} \emph{$\delta$-FDP.}
  Given an estimate of the support $\hat{S}$ obtained from a particular
  inference procedure, the false discovery proportion with parameter
  $\delta > 0$, denoted $\delta$-FDP is the ratio of the number selected
  features that are at a distance more than $\delta$ from any feature of the
  support, divided by the number of selected features:
\begin{equation}\label{eq:delta_fdp}
\delta \text{-FDP} = \frac{\#\{k \in \hat{S}: \forall j \in S, \ d(j,k) > \delta\}}{\#\{k \in \hat{S}\}}
\end{equation}
\end{df}
One can notice that for $\delta = 0$, the FDP and the $\delta$-FDP refer to same
quantity \ie $0\text{\emph{-FDP}} = FDP$.
\begin{df}\label{df:FDR} \emph{False Discovery Rate (FDR) and $\delta$-FDR.}
  The false discovery rate and the false discovery rate with parameter
  $\delta > 0$ which is denoted by $\delta$-FDR are respectively the
  expectations of the FDP and the $\delta$-FDP:
\begin{equation}\label{eq:fdr}
\begin{split}
\text{FDR} & = \mathbb{E}[\text{FDP}] \enspace ,\\
\delta\text{-FDR} & = \mathbb{E}[\delta\text{-FDP}] \enspace . \\
\end{split}
\end{equation}
\end{df}
%
\subsection{Knockoff Inference}
\label{sec:knockoffs}
%
Initially introduced by \cite{barber_controlling_2015} to identify
variables in genomics, the knockoff filter is an FDP control approach for
multivariate models.
This method has been improved to work with mildly high-dimensional settings
in \cite{candes_panning_2016}, leading to the so-called model-X knockoffs:
%
%
\begin{df}
  \label{df:x_knockoff} \emph{Model-X knockoffs \cite{candes_panning_2016}.}
  The model-X knockoffs for the family of random variables
  \( \*X = (X_1, \dots X_p) \) are a new family of random variables
  \( \tilde{\*X} = (\tilde{X}_1, \dots, \tilde{X}_p) \) constructed to satisfy
  the two properties:
\begin{enumerate}
\item For any subset \( \cK \subset \discset{1, \dots, p} \):
    $(\*X, \tilde{\*X})_{\text{swap}(\cK)} \stackrel{d}{=} (\bX, \tilde{\bX})$,\\
  where the vector \( (\*X, \tilde{\*X})_{\text{swap}(\cK)} \) denotes the swap of
  entries \( X_j \) and \( \tilde{X}_j \), \( \forall j \in \cK \)
\item \( \tilde{\bX} \independent \by \mid \bX \) where \( \by \) is the response
  vector.
\end{enumerate}
\end{df}
In a nutshell, knockoff procedure first creates extra null variables that have
a correlation structure similar to that of the original variables.
A test statistic vector is then calculated to measure the strength of the
original versus its knockoff counterpart.
An example of such statistic is the lasso-coefficient difference (LCD) that we
use in this paper:
\begin{df}\label{df:knockoffs_procedure}
\emph{Knockoff procedure with intermediate p-values
\cite{barber_controlling_2015,candes_panning_2016}.}
\begin{enumerate}
\item Construct knockoff variables, produce matrix concatenation:
$[\mathbf{X}, \mathbf{\tilde{X}}] \in \mathbb{R}^{n \times 2p}$.
\item Calculate \textit{LCD} by solving
$$\min_{\mathbf{w} \in \mathbb{R}^{2p}} \dfrac{1}{2} \Vert \mathbf{y} - [\mathbf{X}, \mathbf{\tilde{X}}] \mathbf{w} \Vert^2_2 + \lambda \Vert \mathbf{w} \Vert_1 \enspace , $$
and then, for all $j \in [p]$, take the difference
 $ z_j = \abs{\hat{\mathbf{w}}_j(\lambda)} - \abs{\hat{\mathbf{w}}_{j+p}(\lambda)}$.
%
\item Compute the p-values $p_j$, for $j \in [p]$:
\begin{equation}\label{eq:pval}
  p_j = \dfrac{\# \{k: z_k \leq -z_j \}}{p} \enspace .
\end{equation}
\item Derive q-values by Benjamini-Hochberg procedure: $(q_j)_{j\in [p]} = \text{BHq} \left( (p_j)_{j\in[p]} \right)$
\item Given a desired FDR level $\alpha \in (0,1)$:
 $ \hat{S} = \{j: q_j \leq \alpha \}$.
\end{enumerate}
\end{df}

\begin{rk}\label{rk1}
  The above formulation is distinct from that of \cite{barber_controlling_2015,candes_panning_2016}, but it is equivalent.
  We use it to introduce the intermediate variables $p_j$ for all ${j\in[p]}$.
\end{rk}
%
%
%
%
%
Our first contribution is to extend this procedure computing $q_j$ by
aggregating different draws of knockoffs before applying the Benjamini-Hochberg
(BHq) procedure.
More precisely, we first compute $B$ draws of knockoff variables and, using
\Cref{eq:pval}, we derive the corresponding p-values $p_j^{(b)}$, for all
$j \in [p]$ and $b \in [B]$.
Then, we aggregate them for each $j$
in parallel, using the quantile aggregation procedure introduced in
\cite{meinshausen_p_values_2009}:
\begin{equation}\label{eq:ko_agg}
\forall j \in [p], ~ p_j = \text{quantile-aggregation} (\{ p_j^{(b)} : b \in [B] \})
\end{equation}
We then proceed with the fourth and fifth steps of the knockoff procedure
described in Def. \ref{df:knockoffs_procedure}.
%
%
\subsection{Dimension reduction}
\label{sec:dimension reduction}
%
Knockoff (KO) inference is intractable in high-dimensional settings, as knockoff
generation requires the estimation and inversion of covariance matrices of size
$(2p \times 2p)$.
Hence we leverage data structure by introducing a clustering step that reduces
data dimension before applying KO inference.
As in \cite{chevalier_ensembling_2018}, assuming the features' signals are
spatially smooth, it is relevant to consider a spatially-constrained clustering
algorithm.
By averaging the features with each clustering, we reduce the number of
parameters from $p$ to $q$, the number of clusters, where $q \ll p$.
KO inference on cluster-based signal averages will be referred to as
Clustered Knockoffs (CKO).
However, it is preferable not to fully rely on a particular clustering, as a
small perturbation on the input data has a dramatic impact on the clustering
solution.
We followed the approach used in \cite{hoyos-idrobo_frem_2018} that aggregates
solutions across random clusterings.
More precisely, they build $C$ different clusterings from $C$
different random subsamples of size $\lfloor 0.7n \rfloor$ from the
full sample $\*X$, but always using  the same clustering algorithm.
%
\subsection{The Ensemble of Clustered Knockoff Algorithm}
\label{sec:ECKO_algorithm}
\begin{figure}[t]
    \centering
    \includegraphics[width=0.95\linewidth]{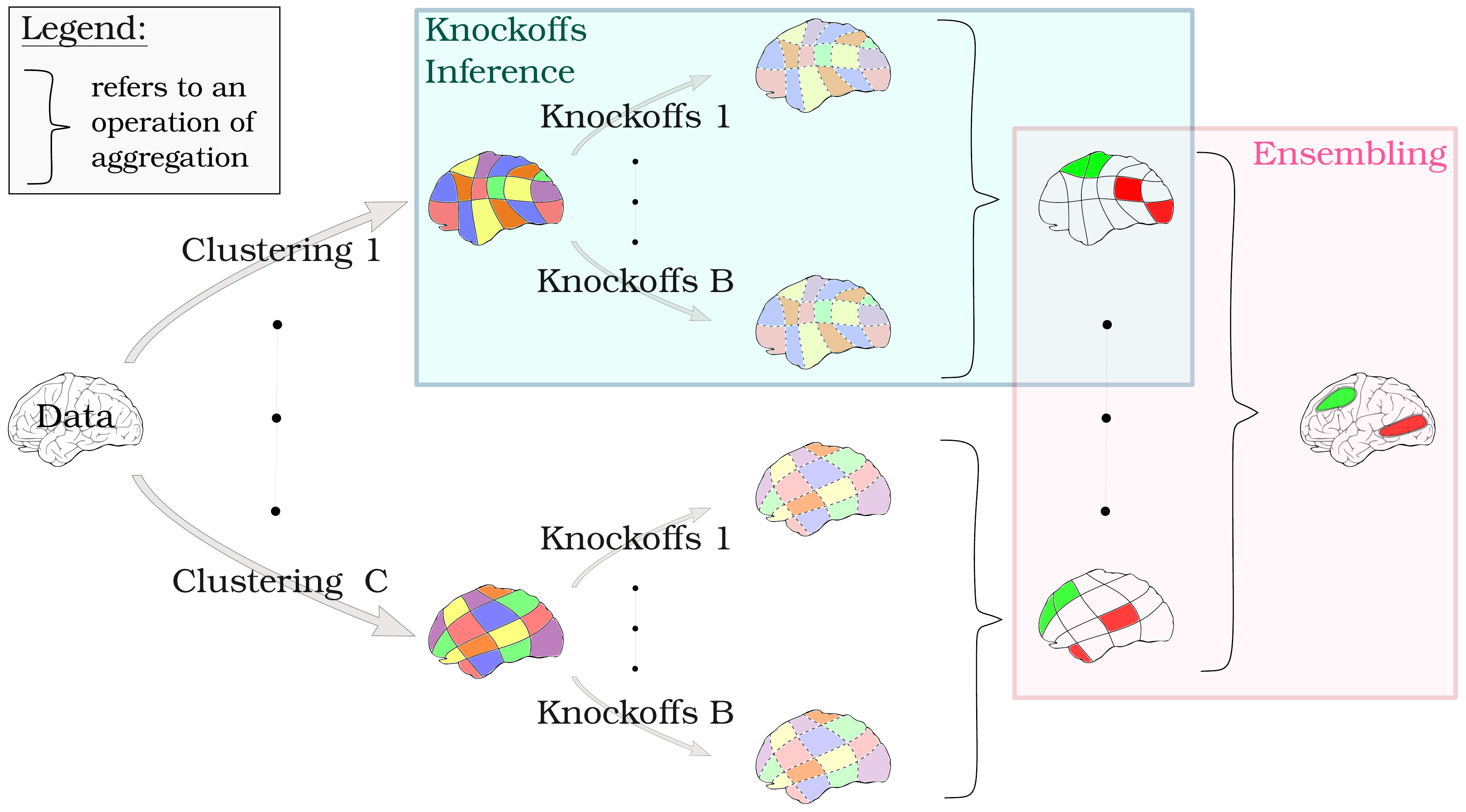}
    \caption{Representation of the ECKO algorithm. To create a stable inference
      result, we introduce ensembling steps both  within each cluster level and at the voxel-level, across clusterings.}
    \label{fig:schema}
\end{figure}

%
The problem is to aggregate the q-values obtained across CKO runs on different
clustering solutions.
To do so, we transfer the q-values from clusters (group of voxels) to features
(voxels): given a clustering solution $c \in [C]$, we assign to each voxel the
q-value of its corresponding cluster.
More formally, if, considering the $c$-th clustering solution, the $k$-th voxel
belongs to the $j$-th cluster denoted by $G^{(c)}_j$ then the q-value
$\tilde{q}^{(c)}_k$ assigned to this voxel is:
$\tilde{q}^{(c)}_k = q_j^{(c)}$ if $k \in G^{(c)}_j$.
This procedure hinges on the observation that the FDR is a
resolution-invariant concept ---it controls the \textit{proportion}
of false discoveries.
In the worst case, this results in a spatial inaccuracy of $\delta$ in the
location of significant activity, $\delta$ being the diameter of the clusters.
Finally, the aggregated q-value $\tilde{q}_k$ of the $k$-th voxel is the
average of the q-values $\tilde{q}^{(c)}_k$, $c \in [C]$, received across
C different clusterings: given the FDR definition \Cref{eq:fdr},
FDPs can naturally be averaged.
The algorithm is summarized in Alg. \ref{alg:ECKO} and
represented graphically in \Cref{fig:schema}.

\begin{algorithm}[h]
\DontPrintSemicolon
\SetNoFillComment
\SetKwInOut{Input}{input}
\SetKwInOut{Output}{output}
\SetKwInOut{Parameter}{param}
\SetAlgoLined
\Input{
  Data matrix \( \bX_{\tt{init}} \in \bbR^{n \times p}\),
  response vector \( \by_{\tt{init}} \in \bbR^{n} \); \\
  Clustering object \(\texttt{Ward}(\cdot) \); \\
}
\Parameter{\( q=500, B=25, C=25,\) \texttt{fdr} - Nominal FDR threshold; \\}
\BlankLine
\For{\( c = 1, \dots, C \)}
{
  \BlankLine
  \( X^{(c)}_{\texttt{init}} = \texttt{resample}(X_{\texttt{init}}) \) \;
  \( X^{(c)}_{\texttt{clustered}} = \texttt{Ward}(q, X^{(c)}_{\texttt{init}}) \) \;
  \BlankLine
  \For{\( b = 1, \dots, B \)}
  { \(\forall j = 1, \dots, q: \) \;
   \( z_j^{(b,c)} \gets \texttt{Knockoffs}( X^{(c)}_{\texttt{clustered}}, \by_{\texttt{init}}, \texttt{fdr})  \) \;
    \( p_j^{(b,c)} \gets \dfrac{\#\{k \in [q]: z_k^{(b,c)} \leq -z_j^{(b,c)}\}}{p} \) \;
  }
  \(\forall j = 1, \dots, q: \) \;
  \( p_j^{(c)} \gets \texttt{Aggregated}(p^{(b, c)}_j, b\in[B]) \) \;
  \( q_j^{(c)} \gets \texttt{BHq\_corrected}(p_j^{(c)}) \) \;
  \( \forall k = 1, \dots, p: \) \;
  \( \tilde{q}^{(c)}_k \gets q_j^{(c)} \ \text{if} \ k \in G^{(c)}_j \) \;
}
  \BlankLine
\( \forall k = 1, \dots, p: \) \;
\( \tilde{q}_k \gets \texttt{Average}(\tilde{q}^{(c)}_k, c \in [C]) \quad \) \;
\BlankLine
\Return \( \hat{S} \gets \discset{k \in [p]: \tilde{q}_k \leq \texttt{fdr}} \)
\caption{Full ECKO algorithm}
\label{alg:ECKO}
\end{algorithm}
%
\subsection{Theoretical Results}
\label{sec:theoretical_results}
%
\paragraph{Ensemble of Clustered Knockoffs (ECKO).}
\begin{theo} \emph{$\delta$-FDR control by the Ensemble of Clustered Knockoffs procedure at the voxel level.}
Assuming that the true model is the one defined in \Cref{linear_model},
using the q-values $\tilde{q}_k$ defined in
\Cref{sec:ECKO_algorithm}, the estimated support
$\hat{S} = \{k : \tilde{q}_k \leq \alpha \}$ ensure that
the $\delta$-FDR is lower than $\alpha$.
\end{theo}
\paragraph{Sketch of the proof}
(details are omitted for the sake of space).  We first establish that the
aggregation procedure yields q-values $q_j^{c}$ that control the FDR.
This follows simply from the argument given in the proof of Theorem
3.1 in \cite{meinshausen_p_values_2009}.
Second, we show that broadcasting the values from clusters ($q$) to voxels
($\tilde{q}$) still controls the FDR, yet with a possible inaccuracy of
$\delta$, where $\delta$ is the supremum of clusters diameters: the $\delta$-FDR
is controlled. This comes from the resolution invariance of FDR and the
definition of $\delta$-FDR.
Third, averaging-based aggregation of the q-values at the voxel level,
controls the $\delta$-FDR. This stems from the definition of the FDR
as an expected value.
\subsection{Alternative approaches}
In the present work, we use two alternatives to the proposed CKO/ECKO approach:
the ensemble of clustered desparsified lasso (ECDL)
\cite{chevalier_ensembling_2018} and the APT framework from \cite{gaonkar_interpreting_2015}.
As we already noted, ECDL is structured as ECKO. The main differences are that
it relies on desparsified lasso rather than knockoff inference and returns
p-values instead of q-values.
The APT approach was proposed to return feature-level p-values for binary
classification problems (though the generalization to regression is
straightforward).
It directly works at the voxel level, yet with two caveats:
\begin{itemize}
\item Statistical control is granted only in the $n \rightarrow \infty$ limit
\item Unlike ECDL and ECKO, it is unclear whether the returned score
  represents marginal or conditional association of the input features with the output.
\end{itemize}
For both ECDL and APT, the returned p-values are converted to q-values
using the standard BHq procedure.
The resulting q-values are questionable, given that BHq is not valid under
negative dependence between the input q-values
\cite{benjamini_controlling_1995}; on the other hand, practitioners rarely check
the hypothesis underlying statistical models. We thus use the procedure in a
black-box mode and check its validity a posteriori.

\section{Experiments}  
\label{sec:experiment}

\subsubsection{Synthetic data.} To demonstrate the improvement of the proposed
algorithm, we first benchmark the method on 3D synthetic data set that resembles
a medical image with compact regions of interest that display some predictive
information.
The size of the weight vector \( \bw \) is \( 50 \times 50 \times 50 \), with 5
regions of interest (ROIs) of size \( 6 \times 6 \times 6 \).
A design matrix \( \bX \) that represents random brain signal is then sampled
according to a multivariate Gaussian distribution. Finally, the response vector
\( \by \) is calculated following linear model assumption with Gaussian noise,
which is configured to have \( SNR \approx 3.6 \), similar to real data
settings.
An average precision-recall curve of 30 simulations is calculated to
show the relative merits of single cluster Knockoffs inference versus ECKO
and ECDL and APT.
Furthermore, we also vary the Signal-to-Noise Ratio (SNR) of the simulation to
investigate the accuracy of FDR control of ECKO with different levels of
difficulty in detecting the signal.

\subsubsection{Real MRI dataset.} We compare single-clustered Knockoffs (CKO),
ECKO and ECDL on different MRI datasets downloaded from the Nilearn library
\cite{abraham_nilearn_2014}. In particular, the following datasets are used:
\begin{itemize}
\item \textbf{Haxby} \cite{haxby_distributed_2001}. In this functional-MRI
  (fMRI) dataset, subjects are presented with images of different objects.
  For the benchmark in our study, we only use the brain signal and responses for
  images related to faces and houses of subject 2 (\(n=216, p=24083 \)).
\item \textbf{Oasis} \cite{marcus_oasis_2007}. The original collection include
  data of gray and white matter density probability maps for 416 subjects aged
  18 to 96, 100 of which have been clinically diagnosed with very mild to
  moderate Alzheimer’s disease. The purpose for our inference task is to find
  regions that predict the age of a subject (\(n=400, p=153809\)).
\end{itemize}
We chose $q=500$ in all experiments for the algorithms that require clustering
step (KO, ECKO and ECDL).
In the two cases, we start with a qualitative comparison of the
returned results.
The brain maps are ternary: all regions outside $\hat{S}$ are zeroed, while
regions in $\hat{S}$ get a value of $+1$ or $-1$, depending on whether the
contribution to the prediction is positive or negative.
For ECKO, a vote is performed to decide whether a voxel is more frequently in a
cluster with positive or negative weight.
%
\section{Results} 
\label{sec:results}
\subsubsection{Synthetic data.} 
\begin{figure}[t!]
  \centering
  \begin{subfigure}[b]{0.48\textwidth}
    \centering
    \includegraphics[width=\textwidth]{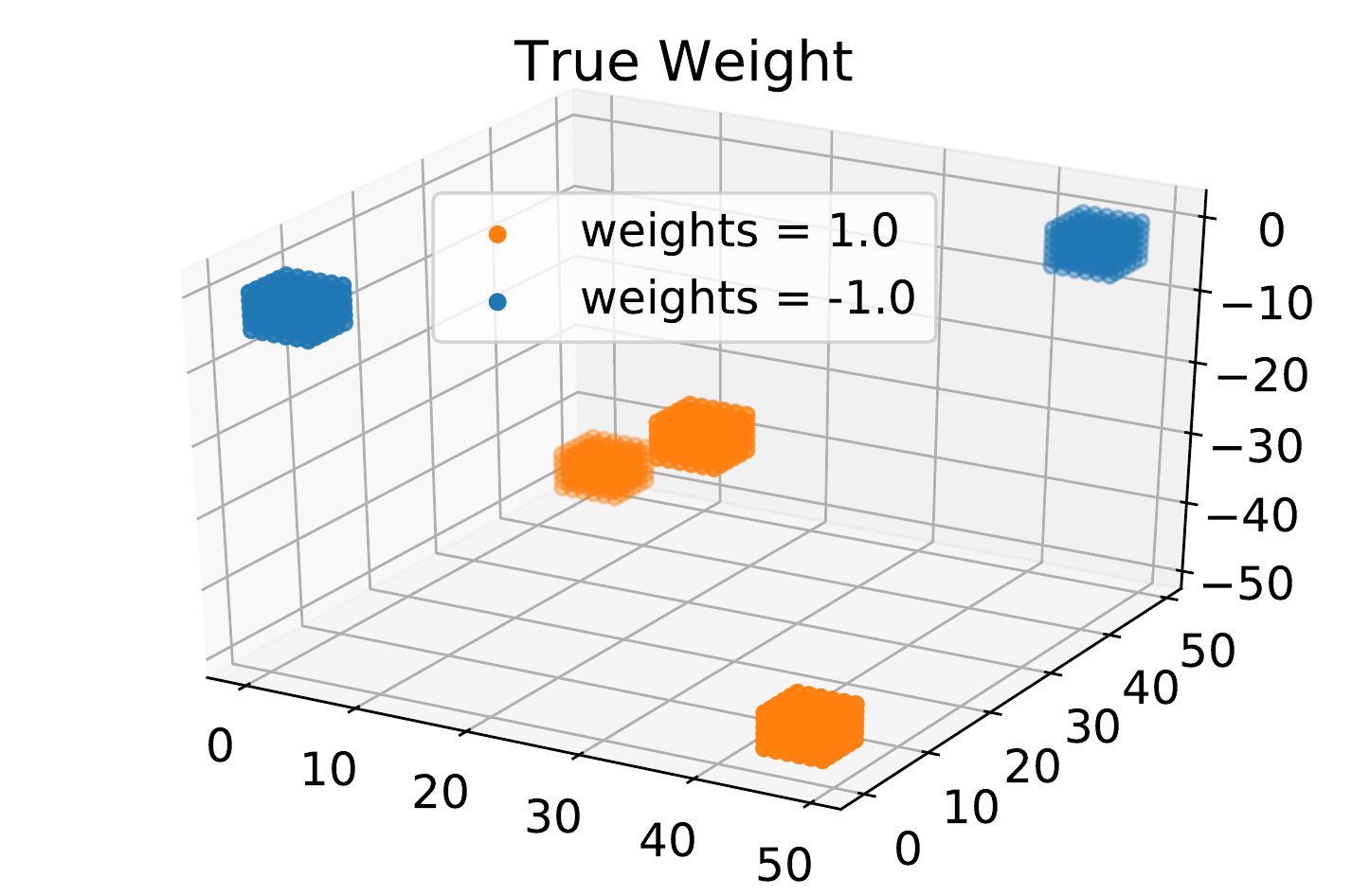}
  \end{subfigure}
  \begin{subfigure}[b]{0.48\textwidth}
    \centering
    \includegraphics[width=\textwidth]{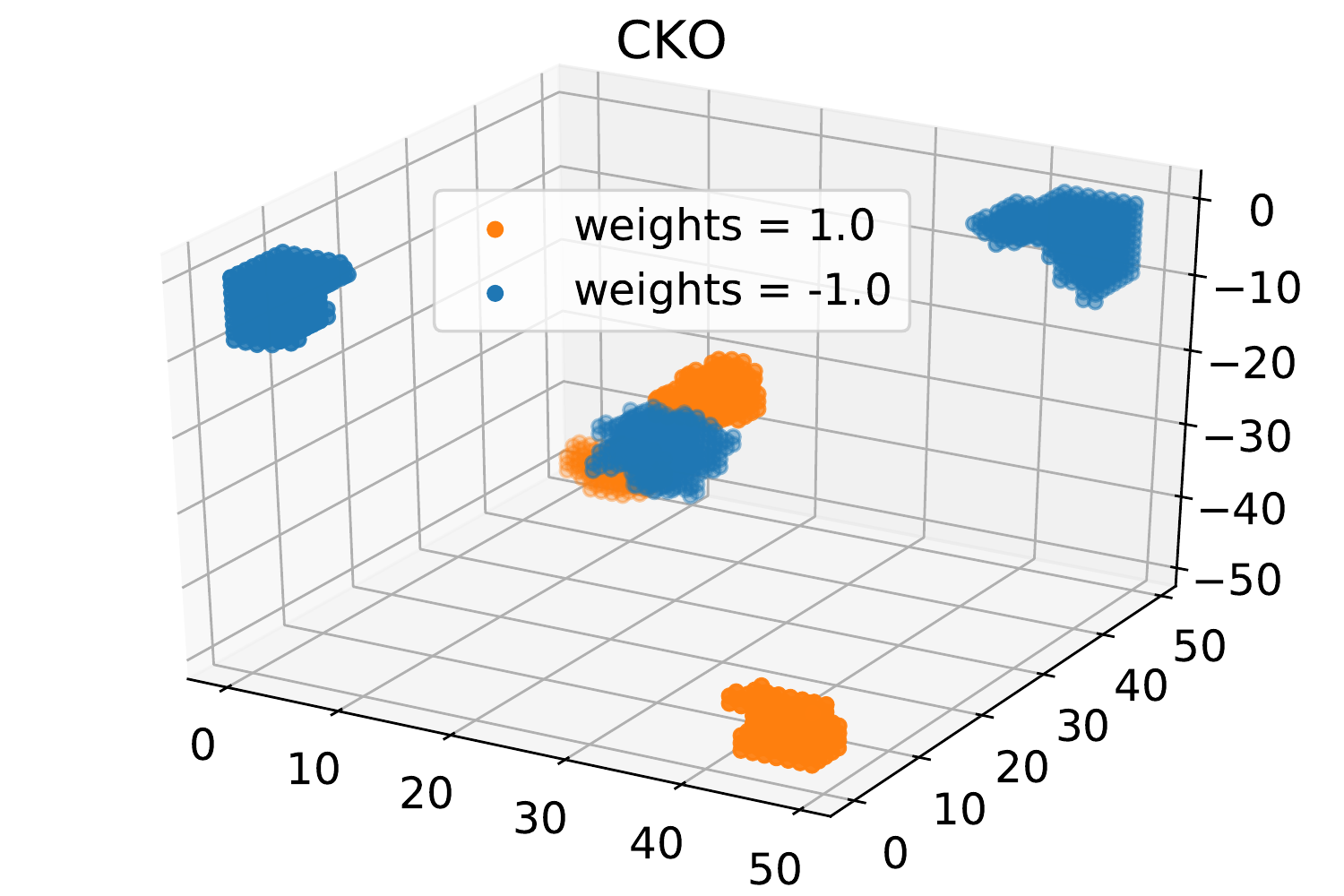}
  \end{subfigure}
  \begin{subfigure}[b]{0.48\textwidth}
    \centering
    \includegraphics[width=\textwidth]{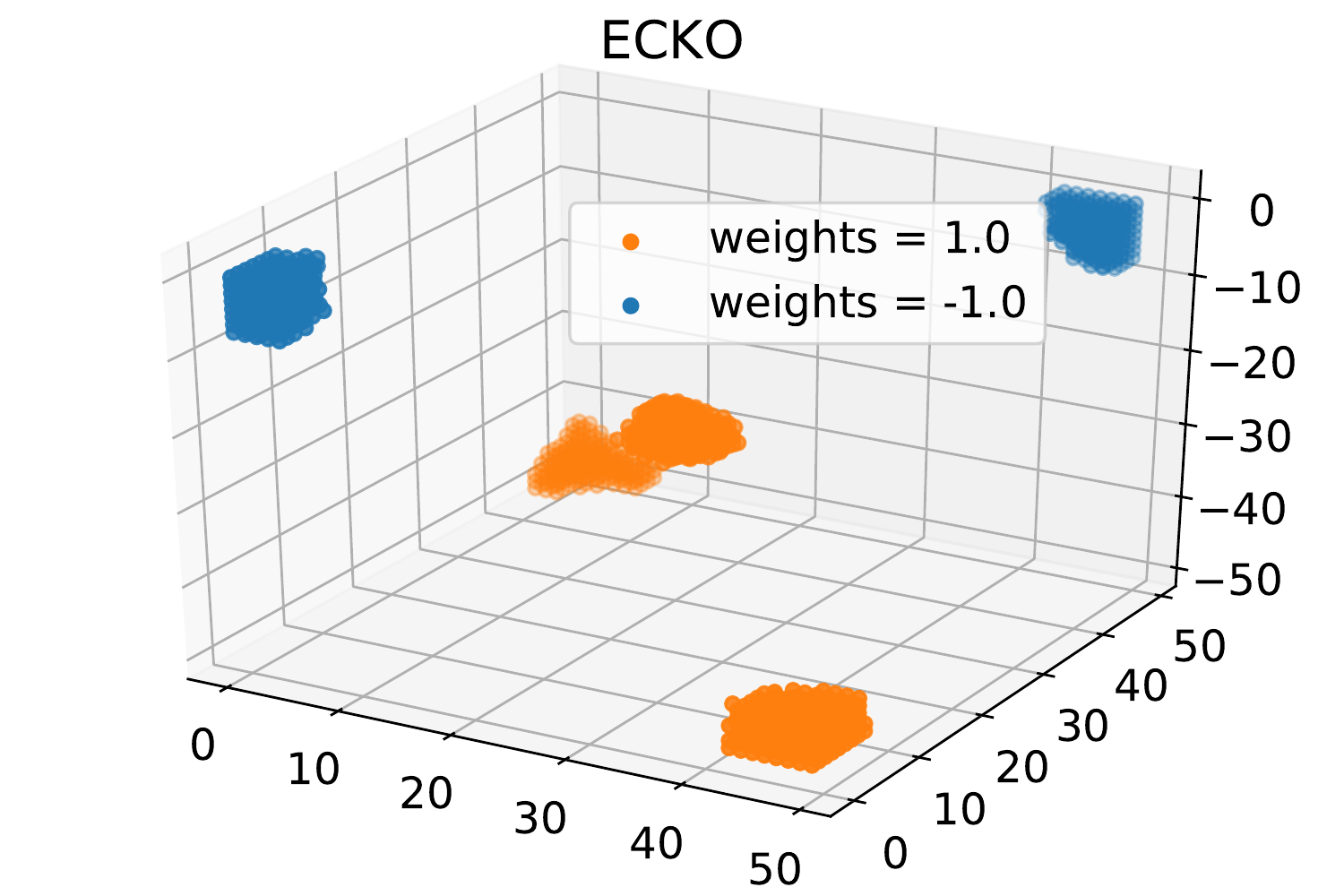}
  \end{subfigure}
  \begin{subfigure}[b]{0.48\textwidth}
    \centering
    \includegraphics[width=\textwidth]{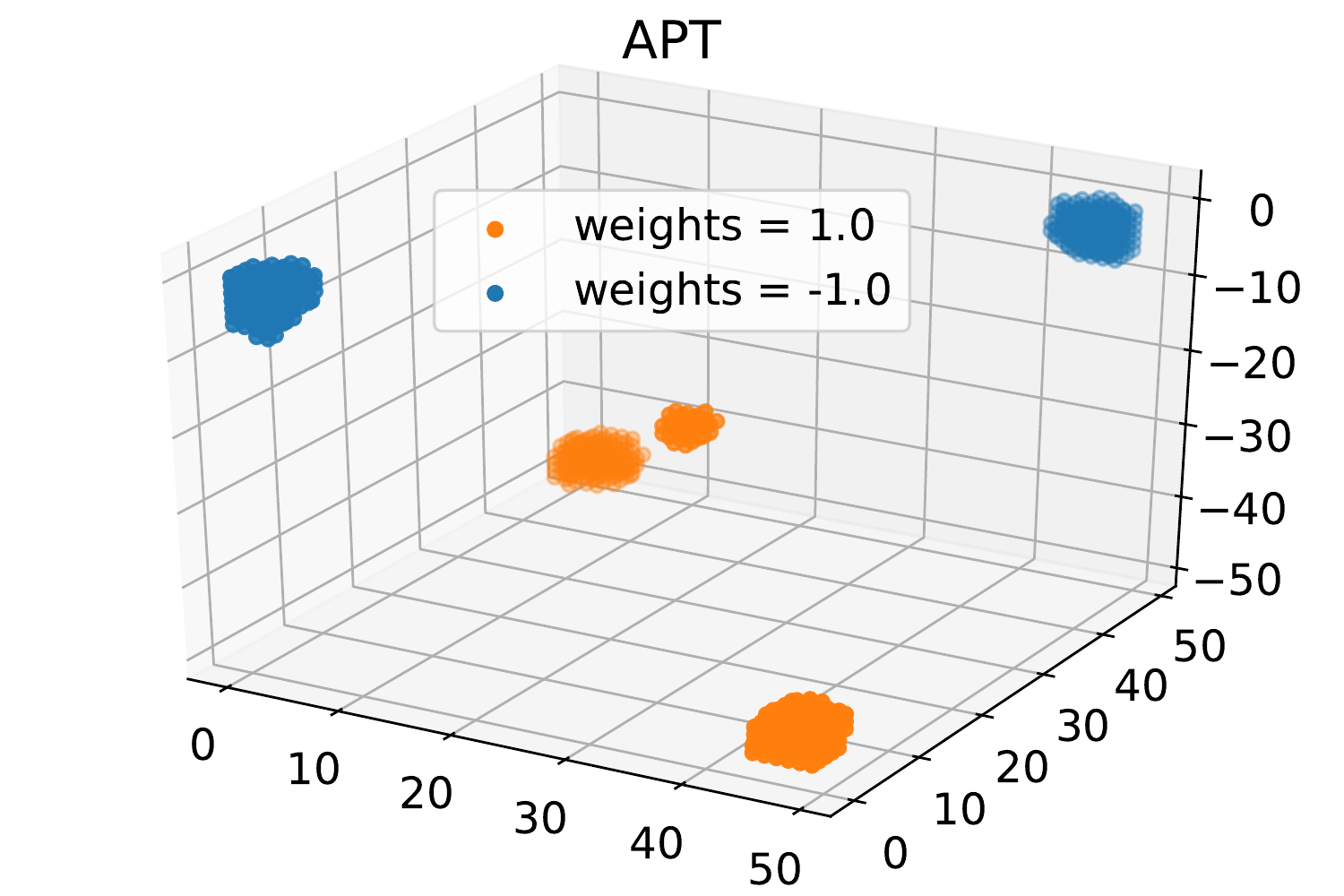}
  \end{subfigure}
  \begin{subfigure}[b]{0.48\textwidth}
    \centering
    \includegraphics[width=\textwidth]{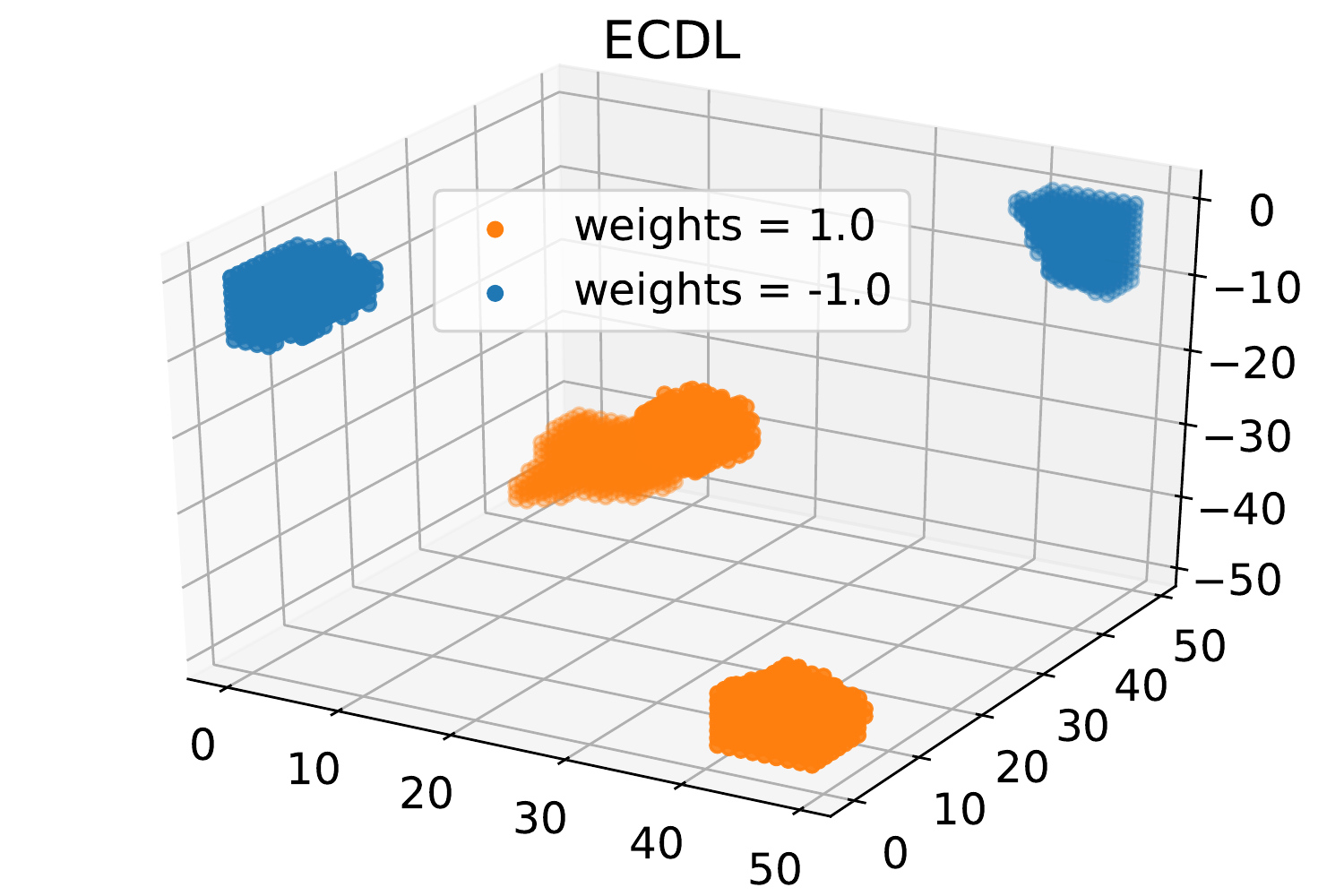}
  \end{subfigure}
  \caption{Experiments on simulated data: Original 3D weight vector (top left)
    and inference results from CKO vs. ECKO. The single CKO run has markedly
    different solutions to the ground truth. Meanwhile, ECKO's solution is
    closer to the ground truth in the sense that altogether, it is more powerful
    than APT and also more precise than ECDL.}
  \label{fig:3D-visualization}
\end{figure}

\begin{figure}[t!]
  \centering
  \begin{subfigure}[b]{0.48\textwidth}
    \centering
    \includegraphics[width=\textwidth]{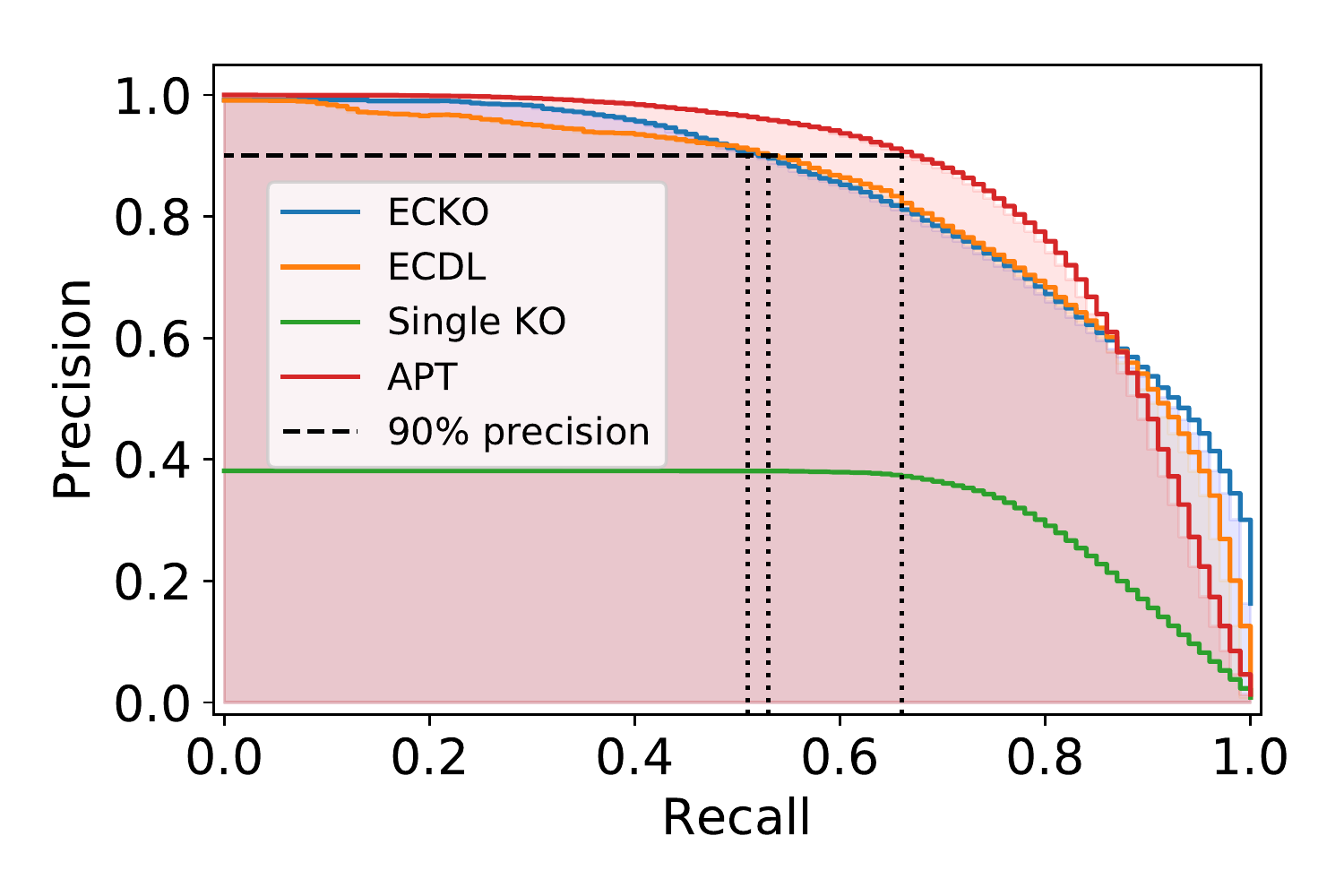}
    \caption{ }
    \label{fig:prec-recall}
  \end{subfigure}
  \begin{subfigure}[b]{0.48\textwidth}
    \centering
    \includegraphics[width=\textwidth]{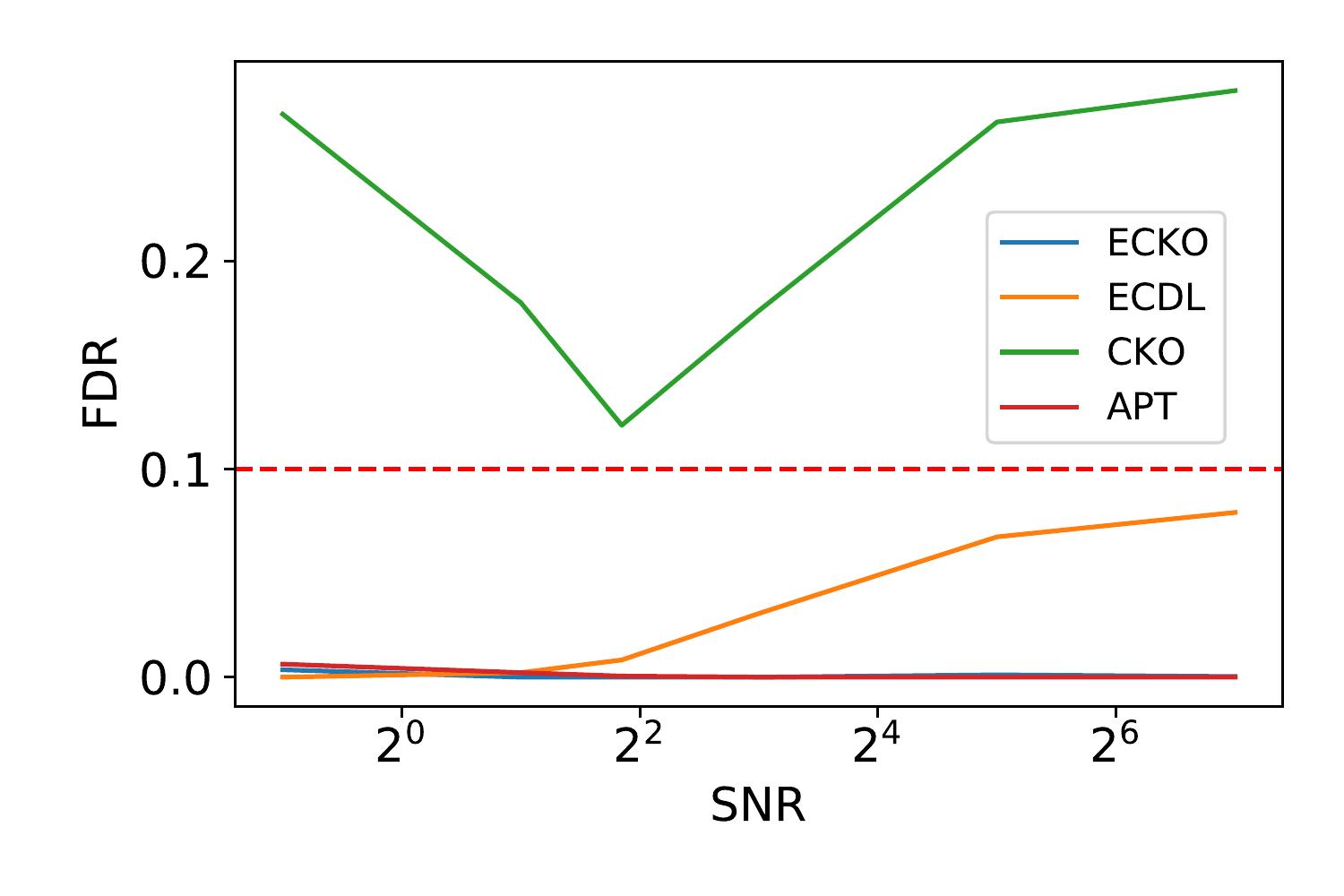}
    \caption{ }
    \label{fig:snr-fdr}
  \end{subfigure}
  \caption{ (a) Average Precision-recall curve (for \(SNR \approx 3.6\)) and (b)
    SNR-FDR curve of 30 synthetic simulations. Nominal FDR control level is
    10\%. ECKO shows substantially better results than CKO and is close to ECDL.
    APT obtains a slightly better Precision-recall curve.
    ECKO, ECDL and APT successfully control FDR under
    nominal level 0.1 where as CKO fails to.}
\end{figure}

A strong demonstration of how ECKO makes an
improvement in stabilizing the single-clustering Knockoffs (CKO) is shown in
Fig. \ref{fig:3D-visualization}.
There is a clear distinction between selection of the orange area at lower right
and the blue area at upper right in the CKO result, compared to the ground
truth.
Moreover, CKO falsely discovers some regions in the middle of the cube.
By Contrast, ECKO's selection is more similar to the true 3D weight
cube.
While it returns a wider selection than ECKO, ECDL also claims more false
discoveries, most visibly in the blue area on upper-left corner.
At the same time, APT returns adequate results, but is more conservative than
ECKO.
Fig. \ref{fig:prec-recall} is the result of averaging 30 simulations for the
3D brain synthetic data.
ECKO and ECDL obtain almost identical precision-recall curve: for a
precision of at least 90\%, both methods have recall rate of around
50\%.
Meanwhile, CKO falls behind, and in fact it cannot reach a precision of over
40\% across all recall rates.
APT yields the best precision-recall compromise, slightly above ECKO and ECDL.
When varying SNR (from \( 2^{-1} \) to \( 2^5 \)) and investigating the average
proportion of false discoveries (\(\delta\)-FDR) made over the average of 30
simulations (Fig. \ref{fig:snr-fdr}), we observe that CKO fails to control
$\delta$-FDR at nominal level 10\% in general.
Note that accurate $\delta$-FDR control would be obtained with larger $\delta$
values, but this makes the whole procedure less useful.
The ECDL controls $\delta$-FDR at low SNR level. However, when the signal is
strong, ECDL might select more false positives.
ECKO, on the other hand, is always reliable ---albeit conservative--- keeping
FDR below the nominal level even when SNR increases to larger magnitude.\\
\\
\textbf{Oasis \& Haxby dataset.}  When decoding the brain signal on
subject 2 of the Haxby dataset using response vector label for watching 'Face
vs. House', there is a clear resemblance of selection results between ECKO and
ECDL.
\begin{figure}[b]
  \centering
  \includegraphics[width=0.49\linewidth]{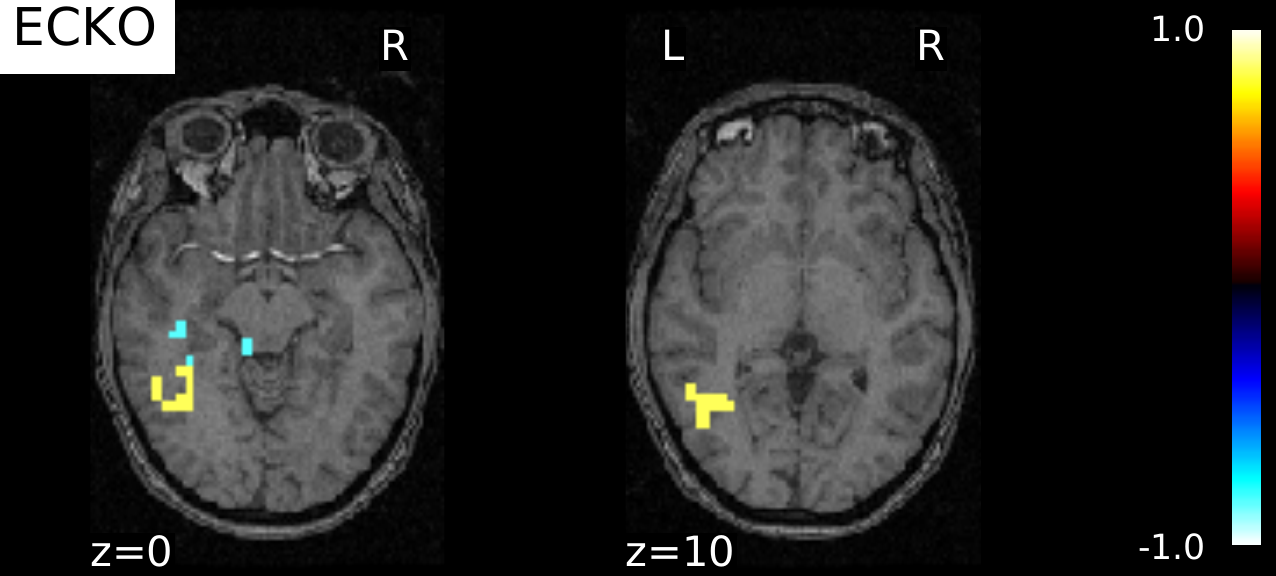}
  \includegraphics[width=0.49\linewidth]{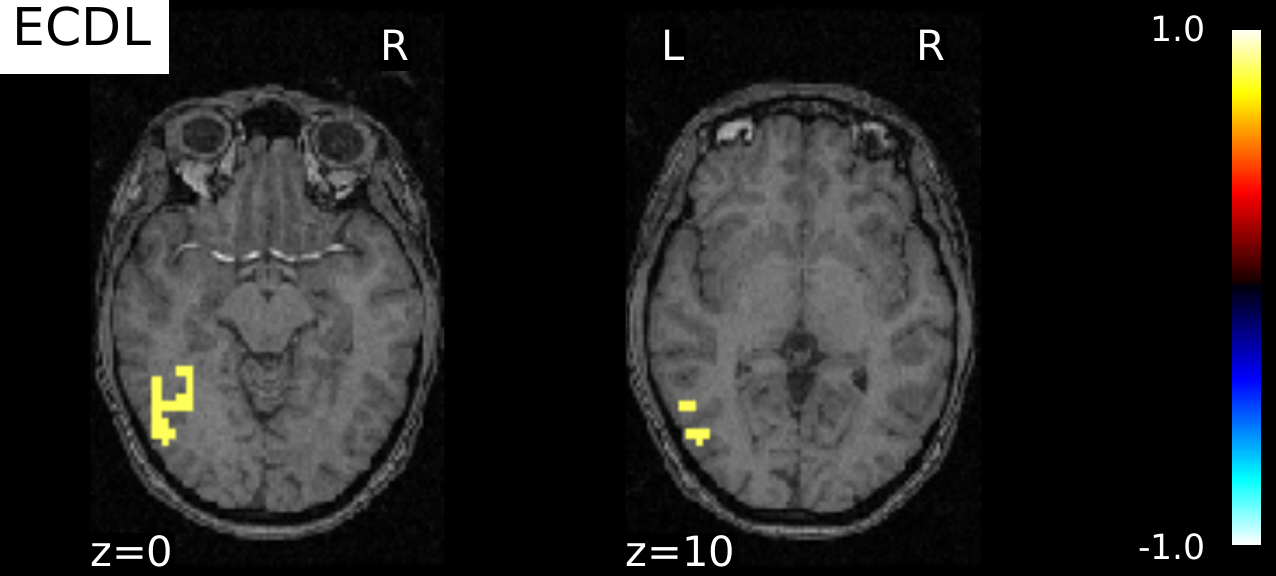}
  \caption{Comparison of results for 2 ensembling clustered inference methods on
    Haxby dataset, nominal FDR=0.1. The results are similar to a large
    extent. No voxel region is detected by APT, therefore we omit to show the
    selection outcome of the method.}
  \label{fig:haxby-ecko-ecdl}
\end{figure}
Using an FDR threshold of 10\%, both algorithms select the same area
(with a difference in size), namely a face responsive region in the
ventral visual cortex, and agree on the sign of the effect.
However, on Oasis dataset, thresholding to control the FDR at 0.1
yields empty selection with ECDL and APT, while ECKO still selects
some voxels.
This potentially means that ECKO is statistically more powerful than
ECDL and APT.
\begin{figure}[t]
  \centering
  \includegraphics[width=0.7\linewidth]{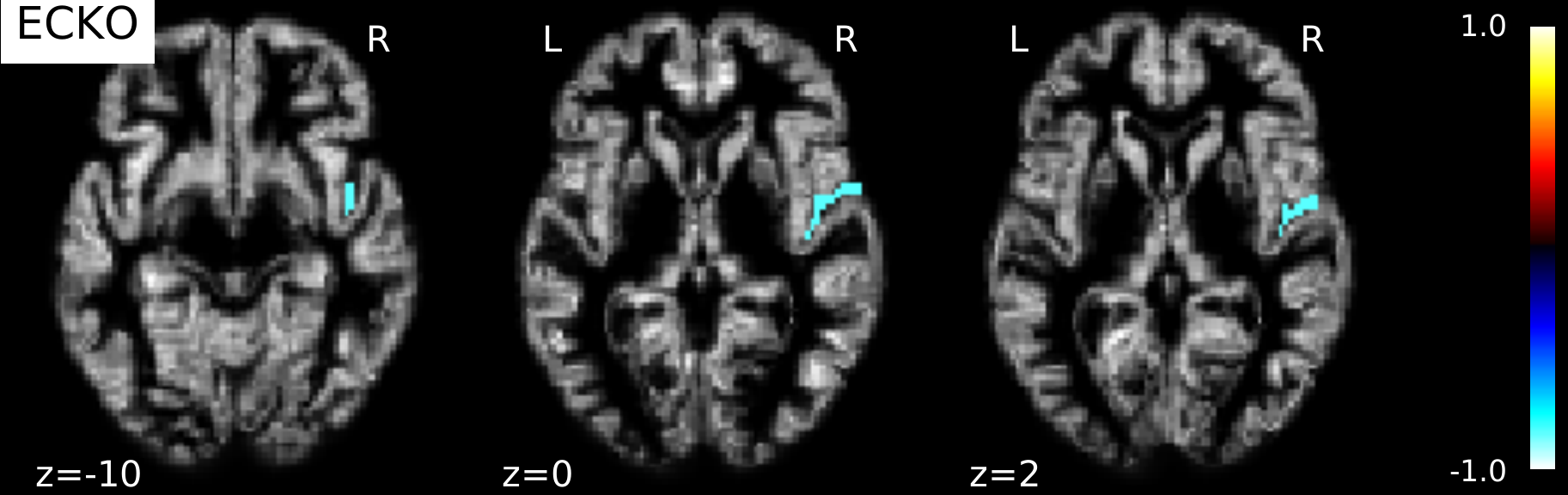}
  \caption{Results of ECKO inference on Oasis dataset, nominal FDR=0.1. ECKO is
    the only method to detect significant regions. The temporal region detected
    by ECKO would be detected by other approaches using a less conservative
    threshold.}
  \label{fig:oasis-ecko}
\end{figure}
%

\section{Conclusion}  
\label{sec:conclusion}

\noindent In this work, we proposed an algorithm that makes False Discovery Rate
(FDR) control possible in high-dimensional statistical inference.
The algorithm is an integration of clustering algorithm for dimension reduction
and aggregation technique to tackle the instability of the original
knockoff procedure.
Evaluating the algorithm on both synthetic and brain imaging datasets
shows a consistent gain of ECKO with respect to CKO in both FDR
control and sensitivity.
Furthermore, empirical results also suggest that the procedure achieves
non-asymptotic statistical guarantees, yet requires the
$\delta$-relaxation for FDR.\\
\indent The number of clusters represents a bias-variance trade-off: increasing
it can reduce the bias (in fact, the value of $\delta$), while reducing it
improves the conditioning for statistical inference, hence the sensitivity of
the knockoff control. We set it to 500 in our
experiments. Learning it from the data is an interesting research direction.\\
\indent We note that an assumption of independence between hypothesis
tests is required for the algorithm to work, which is often not the
case in realistic scenarios.
Note that this is actually the case for all FDR-controlling procedures
that rely on the BHq algorithm.
As a result, making the algorithm work with relaxed assumption is a potential
direction for our future study.
Furthermore, the double-aggregation procedure makes the algorithm more
expensive, although it results in embarrassingly parallel loops.  An
interesting challenge is to reduce the computation cost of this
procedure.
Another avenue to explore for the future is novel generative schemes for
knockoff, based e.g. on deep adversarial approaches. \\ \\
\textbf{Acknowledgement} This research is supported by French ANR (project
FASTBIG ANR-17-CE23-0011) and Labex DigiCosme (project
ANR-11-LABEX-0045-DIGICOSME). The authors would like to thank Sylvain Arlot and
Matthieu Lerasle for fruitful discussions and helpful comments.
%
\bibliographystyle{splncs03}
\bibliography{bibliography}
\end{document}